\documentclass[leqno,12pt]{article}

\usepackage{amsmath, amsthm, amscd, amsfonts, amssymb}

\newtheorem{theorem}{Theorem}[section]
\newtheorem{lemma}[theorem]{Lemma}
\newtheorem{definition}[theorem]{Definition}

\def\vsp{\vspace*{1,5mm}\\ }
\def\barr{\begin{array}}
\def\earr{\end{array}}
\def\({\left(}
\def\){\right)}

\title{Feedback exact controllability of the nonlinear diffusion equation with controller in coefficients} 
\author{Ioan HRINCA}
\date{Department of Mathematics, Otterbein University, Westerville, OH 43081}

\begin{document}
\maketitle

\begin{abstract}
	One designs a feedback controller for semilinear parabolic equations with controller in diffusion coefficient which steers the initial data $y_0$ in a given state $y_1$ in a finite-time.
\end{abstract}

\section{Introduction}\label{s1}
Consider the semilinear parabolic equation
\begin{equation}\label{e1.1}
\barr{ll}
\displaystyle\frac{\partial  y}{\partial  t}-\Delta(uy)+f(x,y)=0&\mbox{in }Q=(0,T)\times{\Omega}\vsp
y(0,x)=y_0(x)&x\in{\Omega}\vsp
y(t,x)=0&(t,x)\in\Sigma=(0,T)\times\partial {\Omega}\earr 
\end{equation} 
where ${\Omega}$ is a bounded and open domain of $\mathbb{R}^d$, $d\ge1$, with sufficiently smooth boundary $\partial {\Omega}$ (of class $C^1$ for instance) and $u=u(t,x)$ is an {\it input controller} which, in order to preserve the parabolic character of the system, is taken in the class
$$\mathcal{U}=\{u\in L^\infty (Q);\ u(t,x)\ge0,\mbox{ a.e. }(t,x_)\in Q\}.$$
Herein, $f:{\Omega}\times\mathbb{R}\to\mathbb{R}$ is a continuous, monotonically nondecreasing function in $y$, and the initial data $y_0$ is taken in the space $L^\infty ({\Omega})$.

The exact controllability problem in time $T$ for \eqref{e1.1} can be briefly\break   described as follows:
\begin{itemize}
	\item[(P)]
{\it Given $y_1\in L^\infty ({\Omega})$, find a controller $u\in\mathcal{U}$ such that
$$y^u(T,x)=y_1(x),\mbox{ a.e. }x\in{\Omega}.$$}
\end{itemize}

\noindent (Herein, $y^u$ is the solution to \eqref{e1.1} in a certain weak sense to be made precise later on.)

Roughly speaking, this is a control problem in coefficients for a semilinear parabolic equation and one should mention that the standard techniques for exact controllability of linear parabolic systems (\cite{6}) are not applicable in the present situation.

Following \cite{3} (see also \cite{2}), we shall treat herein   problem (P) by looking for the controller $u$ in the feedback form
\begin{equation}\label{e1.2}
u=\Phi(x,y)\end{equation}
where $\Phi:{\Omega}\times\mathbb{R}\to\mathbb{R}$ is chosen in such a way that the corresponding closed loop-system
\begin{equation}\label{e1.3}
\barr{ll}
\displaystyle\frac{\partial  y}{\partial  t}-\Delta(\Phi(x,y)y)+f(x,y)=0&\mbox{in }Q\vsp
y(0,x)==y_0(x),&x\in{\Omega}\vsp\
y(t,x)=0,&(t,x)\in\Sigma
\earr
\end{equation}
is well-posed, and its solution $y$ satisfies
$$y(T,x)=y_1(x),\ x\in{\Omega}.$$
Taking into account \cite{2}, \cite{3}, the best candidate for $\Phi$ is the function
\begin{equation}\label{e1.4}
\Phi(x,r)=aH(r-y_1(x)),\ r\ge0,\ x\in{\Omega},\end{equation}
where $a>0$ and $H:\mathbb{R}\to\mathbb{R}$ is the multi-valued Heviside function
\begin{equation}\label{e1.5}
H(r)=\left\{\barr{rll}
1,&\mbox{if}&r>0\\
{[0,1]},&\mbox{if}&r=0\\
0,&\mbox{if}&r<0
\earr\right.\end{equation}
In a few words, the main result of this work is that if $y_0,y_1\in L^\infty ({\Omega})$ and
\begin{equation}\label{e1.6}
y_0(x)\ge y_1(x)\ge0,\ y_0(x)-y_1(x)\ge\alpha>0,\mbox{ a.e. }x\in{\Omega}.
\end{equation}
In other words, $u$ given by \eqref{e1.2} is a real time controller for the exact controllability problem (P). 
Then the closed loop system \eqref{e1.3} has a mild solution $y$ (in a sense to be made precise later on) and $y_0(T,x)=y_1(x)$, a.e. $x\in{\Omega}$. It should be mention, see \cite{2}, \cite{3}, that equation \eqref{e1.3} was used in the case $f\equiv0$ as a model for the {\it self-organized criticality}. Herein, the nonlinear term $f=f(x,y)$ is a potential force, which intervenes on the self-organized criticality process. 

\bigskip\noindent{\bf Notation.} For $1\le p\le\infty $, $L^\infty ({\Omega})$ is the classical space of $p$-integrable functions on $\
ooo$ with the norm denoted $|\cdot|_p$. By $W^{1,p}({\Omega})$ we denote the Sobolev space
$$\left\{u\in W^{1,p}({\Omega});\ \frac{\partial  u}{\partial  x_i}\in L^p({\Omega}),\ i=\overline{1,d}\right\}$$
and
$$W^{1,p}_0({\Omega})=\left\{u\in W^{1,p}({\Omega});\ u\big|_{\partial {\Omega}}=0\right\}.$$
We set 
$$\nabla u=\left\{\frac{\partial  u}{\partial  x_i}\right\}^d_{i=1},\ \Delta u=\sum^d_{i=1}\frac{\partial ^2u}{\partial  x^2_i}.$$
The nonlinear operator (eventually multi-valued) $A:L^1({\Omega})\to L^1({\Omega})$ is said to be accretive if $(I+\lambda  A)^{-1}$ is injective, and Lipschitzian with Lipschitz constant $1$ or, equivalently,  
$$\int_{\Omega}(v_1(x)-v_2(x)){\rm sign}(y_1(x)-y_2(x))dx\ge0,\ \forall  v_i\in A(y_i),\ i=1,2.$$
The operator $A$ is said to be $m$-accretive if it is accretive and $R(I+\lambda  A)=L^1({\Omega}),$ $\forall \lambda >0$ (see \cite{1}). Denote by $C([0,T];L^1({\Omega}))$ the space of continuous functions $u:[0,T]\to L^1({\Omega})$.

\section{The main result}\label{s2}
\setcounter{equation}{0}

We shall assume throughout this work that
\begin{equation}\label{e2.1}
y_0,y_1\in L^\infty ({\Omega})\cap H^1_0(\Omega),\ y_0(x)\ge y_1(x)\ge\alpha,\mbox{ a.e. }x\in\Omega,\end{equation}
where $\alpha$ is a positive constant.

The function $f=f(x,r):{\Omega}\times\mathbb{R}\to\mathbb{R}$ is assumed continuous and monotonically nondecreasing in $r$, measurable in $x\in \Omega$ and
\begin{equation}\label{e2.2}
0\le f(x,r)\le\alpha_0 r,\ \forall  r\ge0\ f(x,0)=0.\end{equation}

\begin{definition}\label{d2.1}\rm The function $y\in C([0,T];L^1(\Omega))$ is said to be a {\it mild solution} to the boundary value problem \eqref{e1.3} if
\begin{equation}\label{e2.3}
y(t)=\lim_{h\to\infty }y_h(t)\mbox{ in }L^1(\Omega)\mbox{ uniformly in }t,\end{equation}
where
\begin{equation}\label{e2.4}
y_h(t)=y^{j+1}_h,\ \forall  t\in[jh,(j+1)h),\ j=0,1,...,N_h=\left[\mbox{$\frac Th$}\right]\end{equation}
\begin{equation}\label{e2.5}
\barr{l}
y^{j+1}_h-h\Delta(\Phi(x,y^{j+1}_h)y^{j+1}_h)+hf(x,y^{j+1}_h)=y^j_h\mbox{ on }\Omega,\vspace*{3mm}\\
y^0_h(x)=y_0(x),\ x\in\Omega.\earr\end{equation}
\end{definition}

Theorem \ref{t2.2} below is the main result.

\begin{theorem}\label{t2.2} Under hypotheses \eqref{e2.1}--\eqref{e2.2} there is a unique mild solution $y$ to \eqref{e1.3}, and it satisfies also
\begin{equation}\label{e2.6}
y(T,x)=y_1(x),\mbox{a.e. }x\in\Omega\end{equation}
if
\begin{eqnarray}
T&=&\eta^{-1}_0\int_\Omega(a^2(y_0(x)-y_1(x))^2+ay_1(x)(y_0(x)-y_1(x))\label{e2.7}\\[2mm]
&&+\int^{y_1(x)}_0f(x,s+y_1(x))ds)dx\nonumber\\[2mm]
\eta_0&=&\gamma^{-2}|y_0-y_1|^{-\frac2{p^*}}_\infty (\alpha a)^{2-\frac2{p^*}}\label{e2.8}\\[4mm]
\gamma&=&\sup\{|v|_{p^*};\|v\|_{W^{1,2}_0}=1\}.\label{e2.9}
\end{eqnarray}
Herein, $p^*=\frac{2d}{d-2},$ if $d>2$, $p^*\in(2,\infty )$ if $d=2$, $p^*=\infty $ if $d=1$. 
\end{theorem}
 
\section{The proof of Theorem \ref{t2.2}}\label{s3}
\setcounter{equation}{0}

It is convenient to reduce equation \eqref{e1.3} by substitution, $z=y-y_1$ to the equation
\begin{equation}\label{e3.1}
\barr{ll}
\displaystyle\frac{\partial  z}{\partial  t}-\Delta(H^*(x,z))+f(x,z+y_1(x))\ni0&\mbox{on }Q\vsp
z(0,x)=z_0(x)=y_0(x)-y_1(x),&x\in\Omega\vsp
z(t,x)=0,&(t,x)\in{\Sigma}.\earr\end{equation}
where
\begin{equation}\label{e3.2}
H^*(x,r)=a(r+y_1(x))H(r),\ x\in\Omega,\ r\ge0.\end{equation}
We consider in the space $L^1(\Omega)$ the operator
\begin{equation}\label{e3.3}
\barr{lcl}
A(z)&=&\{-a\Delta v;v(x)\in H^*(x,z(x))\},\mbox{ a.e. }x\in\Omega\vsp
D(A)&=&\{z\in L^1(\Omega);\Delta v\in L^1(\Omega)\}\earr\end{equation}
where $\Delta$ is taken in the sense of distributions on $\Omega$. Let $B:L^2({\Omega})\to L^2({\Omega})$ be defined by
\begin{equation}\label{e3.4}
B(z)=f(x,z+y_1(x)),\mbox{ a.e. }x\in\Omega,\ z\in L^1(\Omega).\end{equation}
We have (see Lemma 4.1 in \cite{3}).

\begin{lemma}\label{l3.1} The operator $A$ is accretive in $L^1(\Omega)$. Moreover,
\begin{equation}\label{e3.5}
\mathcal{C}=\{z\in L^1(\Omega);z(x)\ge0,\mbox{ a.e. }x\in\Omega\}\subset R(I+\lambda  A),\ \forall \lambda >0.\end{equation}
\end{lemma}
We note also that $B$ is continuous and accretive in $L^1(\Omega)$. Moreover, as easily seen, we have
\begin{equation}\label{e3.6}
\mathcal{C}\subset R(I+\lambda  B),(I+\lambda  B)^{-1}(\mathcal{C})\subset\mathcal{C},\ \forall  \lambda >0.\end{equation}

\begin{lemma}\label{l3.2} We have
\begin{equation}\label{e3.7}
\mathcal{C}\subset R(I+\lambda (A+B)),\ \forall \lambda >0.\end{equation}
\end{lemma}

\begin{proof} One must prove that, for each $g\in\mathcal{C}$ and $\lambda >0$, the equation
\begin{equation}\label{e3.8}
v+\lambda  A(v)+\lambda  B(v)=g\end{equation}
has a solution $v\in D(A)\cap\mathcal{C}$. To this end, we approximate \eqref{e3.8} by
\begin{equation}\label{e3.9}
v_n+\lambda  A(v_n)+\lambda  B_n(v_n)=g,\end{equation}
where $B_n=n\(I-\(I+\frac1nB\)^{-1}\)$, $n\in\mathbb{N}$, is the Yosida approximation of $B$ (see \cite{1}). We may rewrite equation \eqref{e3.9} as
$$(n\lambda +1)v_n+\lambda  A(v_n)=g+n\lambda \(I+\frac1n\,B\)^{-1}(v_n).$$
Equivalently,
\begin{equation}\label{e3.10}
v_n=\(I+\frac\lambda {n\lambda +1}\,A\)^{-1}
\(\frac1{n\lambda +1}\,g+\frac{n\lambda }{n\lambda +1}\(I+\frac1n\,B\)(v_n)\).\end{equation}

Taking into account that, by Lemma \ref{l3.1},
$\(I+\frac\lambda {n\lambda +1}\,A\)^{-1}(\mathcal{C})\subset\mathcal{C}$ and\break $ \(I+\frac\lambda {n\lambda +1}\,A\)^{-1}$ is an $L^1$  contraction on $\mathcal{C}$, it follows by Banach's fixed point theorem that \eqref{e3.10} (and, therefore, \eqref{e3.9}) has a unique solution $v_n\in D(A)\cap\mathcal{C}$. Moreover, taking into account that 
$$\int_{\Omega} B_n(v){\rm sign}\,v\,dx\ge0,\ \int_{\Omega} A(v_n){\rm sign}\,v\,dx\ge0,$$
we get by \eqref{e3.9} that
$$|v_n|_1+|A(v_n)|_1+|B_n(v_n)|_1\le c,\ \forall  n,$$
which implies, in particular, that
$$\sup_n\int_{\Omega}|\Delta v_n(x)|dx<\infty ,$$
and this yields
\begin{equation}\label{e3.11}
\sup_n\|v_n\|_{W^{1,q}_0({\Omega})}<9,\end{equation}
for $q\in\(1,\frac d{d-1}\)$ (see, e.g., \cite{1}, p.110). 
Since $W^{1,q}_0({\Omega})$ is compactly embedded in $L^1{\Omega})$, we infer by  \eqref{e3.11} that $\{v_n\}$ is compact in $L^1({\Omega})$ and so, on a subsequence we have
$$v_n\to v\mbox{ strongly in }L^1\Omega).$$
Then, by \eqref{e3.9} we have
$$A(v_n)\to A(v)=\Delta v\mbox{ in }\mathcal{D}'(\Omega),$$
where$\mathcal{D}'(\Omega)$ is the space of distributions on $\Omega$. Since, by the continuity of $B$ we have also for $n\to\infty $
$$B_n(v_n)\to B(v)\mbox{ in }L^1\Omega),$$
it follows that $v$ satisfies\eqref{e3.8}. It is also clear (by letting $n\to\infty $ in \eqref{e3.9})  that $v\in\mathcal{C}$.\end{proof}

 By Lemma \ref{l3.2} it follows via the Crandall-Liggett theorem (see \cite{1}, p.146) that the Cauchy problem 
\begin{equation}\label{e3.12}
\barr{l}
\displaystyle\frac{dz}{dt}+A(z)+B(z)\ni0,\ t\in[0,T)\vsp 
z(0)=z_0\earr\end{equation}
has, for each $z\in\mathcal{C}$, a unique mild solution $z\in C([0,T];L^1({\Omega}))$, $z(t)\in\mathcal{C}$, $\forall  t\in[0,T]$. 
In other words,
\begin{eqnarray}
&&z(t)=\lim_{h\to0}z_h(t)\mbox{ in }L^1(\Omega),\ \forall  t\in[0,T]\label{e3.13}\\[2mm]
&&z_h(t)=z^{i+1}_h,\ \forall  t\in[ih,(i+1)h),\ i=0,1,...,N_h=\left[\mbox{$\frac Th$}\right]\label{e3.14}\\[2mm]
 &&z^{i+1}_h+hA(z^{i+1}_h)+hB(z^{i+1}_h)\ni z^i_h,\ i=0,1,...,N_h,\label{e3.15}\\[2mm]
&&z^0_h=z_0.\label{e3.16}
\end{eqnarray}
We note that $y=z+y_1$ is just a mild solution to \eqref{e1.3} in the sense of Definition \ref{d2.1}. We pause briefly to note that, since $z_0\in L^\infty (\Omega)$, it follows by \eqref{e3.15} that $z^{i+1}_h\in L^\infty (Q)$,
\begin{equation}\label{e3.17}
|z^{i+1}_h(t)|_\infty \le|z_0|_\infty ,\ \forall  i=0,1,...,N_h,\end{equation}
and, therefore,
\begin{equation}\label{e3.18}
|z_h(t)|_\infty \le|z_0|_\infty ,\ \forall  t\in(0,T).\end{equation}
The argument is standard and we omit it (see \cite{3}). Moreover, by the elliptic regularity, it follows by \eqref{e3.15} that $z^{i+1}_h\in W^{1,2}_0(\Omega)$ and so 
$$z_h(t)\in W^{1,2}_0(\Omega),\ \forall  t\in(0,T).$$
Next, we consider the potential $j=j(x,r)$ of the function $H^*(x,v)$, that is (see \eqref{e3.4})
$$H^*(x,r)=\partial _rj(x,r),\ r\in\mathbb{R},$$
where
$$j(x,r)=\left\{\barr{ll}
\displaystyle\frac a2\,r^2+ary_1(x),&\mbox{if }r\ge0\vsp
0,&\mbox{if }r<0\earr\right.$$
and $\partial _rj(x,t)$ is the subdifferential of $r\to j(x,r)$, that is,
\begin{equation}\label{e3.19}
j(x,r)\le j(x,\bar r)+H^*(x,r)(r-\bar r),\ \forall  r,\bar r\in\mathbb{R}.\end{equation}
We set also
$$g(x,r)=\int^r_0f(x,s+y_1(x))ds,\ \forall  r\ge0,$$
which yields
$$g(x,r)\le g(x,\bar r)+f(x,r+y_1(x))(r-\bar r),\ \forall  r,\bar r\ge0.$$
Then, by \eqref{e3.15} we get
$$\barr{lcl}
\displaystyle\int_{\Omega} j(x,z^h_{i+1}(x))dx&\le&\displaystyle\int_\Omega H^*(x,z^h_{i+1}(x))(z^h_{i+1}(x)-z^h_i(x))dx\vsp
&+&\displaystyle\int_\Omega j(x,z^h_i(x))dx,\earr$$
this yields (by \eqref{e3.15})
$$\barr{lcl}
\displaystyle\int_{\Omega} j(x,z^h_{i+1}(x))dx&\le&\displaystyle\int_\Omega j(x,z^h_i(x))dx
+h\int_\Omega B(z^h_{i+1})(z^h_{i+1}-z^h_i)dx\vsp
&&-ah\displaystyle\int_\Omega|\nabla_xH^*(x,z^h_{i+1}(x))|^2dx,\earr$$
and, therefore, 
$$\barr{lcl}
\displaystyle\int_{\Omega}(j(x,z^h_{i+1}(x))+g(x,z^h_{i+1}(x)))dx
+ah\displaystyle\int_\Omega|\nabla_xH^*(x,z^h_{i+1}(x))|^2dx\vsp
\qquad\le\displaystyle\int_\Omega(j(x,z^h_1(x))+g(x,z^h_1(x)))dx, i=0,1,...,N_h.\earr$$
We set $\eta^h_{i+1}(x)=H^*(x,z^h_{i+1}(x))$ and get
$$\barr{lcl}
\displaystyle\int_{\Omega}(j(x,z^h_{i+1}(x))+g(x,z^h_{i+1}(x)))dx+ah\int_\Omega|\nabla\eta^h_{i+1}(x)|^2dx\vsp
\qquad\le\displaystyle\int_\Omega j(x,z^h_i(x))+g(x,z^h_i(x))dx.\earr$$
This yields
\begin{equation}\label{e3.20}
\barr{l}
\displaystyle\int_{\Omega} j(x,z^h_{i+1}(x))+g(x,z^h_{i+1}(x))dx
+a\ell\displaystyle\sum^{i}_{k=\ell}\int_{\Omega}|\nabla\eta^h_{k+1}(x)|^2dx\vsp
\qquad\le\displaystyle\int_{\Omega}(j(x,z^h_{\ell+1}(x))+g(x_1,z^h_{\ell+1}(x)))dx,\ 0\le\ell\le N_h,\earr\end{equation}
and, therefore,	
\begin{equation}\label{e3.21}
\barr{l}
\displaystyle\int_{\Omega} (j(x,z^h(t,x))+g(x,z^h(t,x)))dx
+a\displaystyle\int^t_s\int_{\Omega} |\nabla\eta^h_{\tau,x}(x)|^2dxd\tau\vsp
\qquad\le\displaystyle\int_{\Omega} j(x,z^h_{\ell+1}(x))+g(x_1,z^h_{\ell+1}(x))dx,\mbox{ for }0\le s\le t<T.\earr\end{equation}
Herein, $\eta^h(\tau,x)=\eta^h_{\ell+1}(x)$ for $\tau=[\ell h,(\ell+1)h]$.

On the other hand, we have by the 
Sobolev--Gagliardo--Nirenberg  theo\-rem  (see \cite{4}, p.278)
$$\(\int_\Omega|\eta^h_{\ell+1}(x)|^{p^*}dx\)^{\frac1{p^*}}
\le\gamma\(\int_\Omega|\nabla\eta^h_{\ell+1}(x)|^2dx\)^{\frac12},$$
where $p^*=\frac{2d}{d-2}$  if $d>2$, $p^*\in(2,\infty )$ if $d=2$, and $p^*=\infty $ if $d=1$, and 
\begin{equation}\label{e3.22}
\gamma=\sup\{|v|_{p^*};\|v\|_{W^{1,2}_0(\Omega)}=1\}.\end{equation}
Then, by \eqref{e3.21} we get
\begin{equation}\label{e3.23}
\hspace*{-2mm}\barr{l}
\displaystyle\int_{\Omega} j(x,z^h(t,x))+g(x,z^h(t,x))dx
+a\gamma^{-1}\!\displaystyle\int^t_s\!\!\(\int_{\Omega}|\eta^h(\tau,x)|^{p^*}dx\)^{\frac2{p^*}}\!\!d\tau\vsp 
\qquad\le\displaystyle\int_{\Omega}(j(x,z^h(s,x))+g(x,z^h(s,x)))dx. \earr\hspace*{-3mm}
\end{equation}
Now, taking into account that
$$\(\int_\Omega|\eta^h(z,x)|^{p^*}dx\)^{\frac2{p^*}}=a^2\(\int_{(z^h>0)}|z^h(z,x))+y_1(x)|^{p^*}dx\)^{\frac2{p^*}},$$
we get by \eqref{e3.23}
$$\barr{l}
\displaystyle\int_{\Omega} (j(x,z^h(t,x))+g(x,z^h(t,x)))dx\vsp
\qquad+a^2\gamma^{-2}\displaystyle\int^t_0\(\int_{[z^h>0]}|z^h(\tau,x)+y_1(x)|^{p^*}dx\)^{\frac2{p^*}}d\tau,\earr$$
and this yields
\begin{equation}\label{e3.24}
\barr{l}
\displaystyle\int_{\Omega} (j(x,z^h(t,x))+g(x,z^h(t,x)))dx\vsp
\qquad+a^2\gamma^{-2}\alpha^2|z_0|^{-\frac2{p^*}}_\infty 
\displaystyle\int^t_s\int_{\Omega} z^h(t,x)dx\vsp
\qquad\le\displaystyle\int_{\Omega} j(x,z^h(s,x))+g(x,z^h(s,x))dx. 
\earr\end{equation}
Let $\psi(t)=\int_\Omega(j(x,z(t,x))+g(x,z(t,x)))dx.$
Then, by \eqref{e3.24} we see that
$$\barr{c}
\psi'(t)+a^2\alpha^2\gamma^{-2}|y_0-y_1|^{-\frac2{p^*}}_\infty 
\(\displaystyle\int_\Omega z(t,x)dx\)^{\frac2{p^*}}\le0,\mbox{ a.e. }t\in(0,T),\vsp 
\psi(0)=\displaystyle\int_\Omega(j(x,z_0(x))+g(x,z_0(x)))dx.\earr$$
Taking into account that $j(t,z)\ge\alpha az$, we get
$$\int_\Omega z(t,x)dx\ge\frac{\psi(x)}{a\alpha},\ \forall  t\in[0,T],$$
and, therefore,
$$\psi'(t)+\eta_0(\psi(t)))^{\frac2{p^*}}\le0,\mbox{ a.e. }t\in(0,T),$$
where
$$\eta_0=\gamma^{-2}(\alpha a)^{2-\frac2{p^*}}|y_0-y_1|^{-\frac2{p^*}}_\infty .$$
This implies that
$$\psi(T)=0\mbox{ if }T=\eta^{-1}_0\int_\Omega j(x,(y-y_0)(x))+g(x,(y_0-y_1)(x))dx.$$
We have shown, therefore, that, for such $T$,
$$z(T,x)=0,\mbox{ a.e. }x\in\Omega,$$and, recalling that $y=z+y_1$, we see that
\eqref{e2.6}--\eqref{e2.7} hold.\hfill$\Box$ 

\bigskip\noindent{\bf Acknowledgements.} This work was done during the visit of the author at Octav Mayer Institute of Mathematics of Romanian Academy.

\end{document}